\documentclass[12pt]{amsart}
\usepackage{amsmath,amsthm,amsfonts,amssymb,mathrsfs}
\date{\today}

\usepackage{hyperref}

\newtheorem{theorem}{Theorem}%[section]

\newtheorem{proposition}[theorem]{Proposition}%[section]
\newtheorem{corollary}[theorem]{Corollary}%[section]

\newtheorem{lemma}[theorem]{Lemma}%[section]
\theoremstyle{definition}
\newtheorem{definition}[theorem]{Definition}

\newtheorem{example}[theorem]{Example}%[section]
\newcommand{\J}{{\mathscr I}}
\newcommand{\D}{{\mathcal D}}

\begin{document}

%\title[On topological inverse semigroup of finite transformations
%of an infinite set]
%{On topological inverse semigroup of finite
%transformations of an infinite set}

\title[Semigroup Closures of Finite Rank Symmetric Inverse Semigroups]
{Semigroup Closures of Finite Rank Symmetric Inverse Semigroups}

\author{Oleg~Gutik}
\address{Department of Mechanics and Mathematics, Ivan Franko Lviv
National University, Universytetska 1, Lviv, 79000, Ukraine}
\email{o\_\,gutik@franko.lviv.ua}

\author{Jimmie~Lawson}
\address{Department of Mathematics,
Louisiana State University,
Baton Rouge, LA 70803, USA}
\email{lawson@math.lsu.edu}

\author{Du\v{s}an~Repov\v{s}}
\address{Institute of Mathematics, Physics and Mechanics, and Faculty of Mathematics and Physics,
University of Ljubljana, Jadranska 19, Ljubljana, 1000, Slovenia}
\email{dusan.repovs@guest.arnes.si}



\begin{abstract}
We introduce the notion of semigroup with a tight
ideal series and investigate their closures in semitopological semigroups,
particularly inverse semigroups with continuous inversion. As a corollary we
show that the symmetric inverse semigroup of finite transformations
$\mathscr{I}_\lambda^n$ of the rank $\leqslant n$ is algebraically
closed in the class of (semi)topological inverse semigroups with
continuous inversion. We also derive related results about
the nonexistence of (partial) compactifications of classes of
semigroups that we consider.
\end{abstract}

\maketitle

%\tableofcontents

\section*{Introduction}

A partial one-to-one transformation on a set $X$ is a one-to-one
function with domain and range subsets of $X$ (including the empty
transformation with empty domain). There is a natural associative
operation of composition on these transformations, $ab(x)=a(b(x))$
wherever defined,  and the resulting semigroup is called the
\emph{symmetric inverse semigroup} $\J_X$ \cite{CP}. The symmetric
inverse semigroup was introduced by V.~V.~Wagner~\cite{Wagner1952}
and it plays a major role in the theory of semigroups. If the
domain is finite and has cardinality $n$, which is then also the
cardinality of the range, the transformation is said to be of
\emph{rank} $n$.  For each $n\geq 0$, the members of $\J_X$ of
rank less than or equal to $n$ form an ideal of $\J_X$, denoted
$\J_X^n$.  (Recall that a nonempty subset $I$ of a semigroup $S$
is an \emph{ideal} if $SI\cup IS\subseteq I$.)  If $X$ and $Y$
have the same cardinality, then these respective semigroups are
isomorphic, and thus we may restrict our attention to a canonical
one from that class, which we take to be the one arising by taking
$X$ to be the cardinal $\lambda$.  We thus label the corresponding
semigroups as $\J_\lambda$ and $\J_\lambda^n$ respectively. The
semigroups $\J_\lambda^n$ for $\lambda$ infinite form a motivating
example for the considerations and developments of this paper.

Many topologists have studied topological properties of
topological spaces of partial continuous maps $\mathscr{PC}(X,Y)$
from a topological space $X$ into a topological space $Y$ with
various topologies such as the Vietoris topology, generalized
compact-open topology, graph topology, $\tau$-topology, and others
(see \cite{Abd-AllahBrown1980, BoothBrown1978,
DiConcilioNaimpally1998, Filippov1996, Hola1998, Hola1999,
KumziShapiro1997, Kuratowski1955}). Since the set of all partial
continuous self-transformations $\mathscr{PCT}(X)$ of the space
$X$ with the operation composition is a semigroup,  many semigroup
theorists have considered the semigroup of continuous
transformations (see surveys \cite{Magill1975} and
\cite{GluskinScheinSnepermanYaroker1977}), or the semigroup of
partial homeomorphisms of an arbitrary topological space (see
\cite{Baird1977, Baird1977a, Baird1977b, Baird1979, Gluskin1977,
Mendes-GoncalvesSullivan2006, Sneperman1962, Yaroker1972}).
Be\u{\i}da~\cite{Beida1980}, Orlov~\cite{Orlov1974, Orlov1974a},
and Subbiah~\cite{Subbiah1987} have considered  semigroup and
inverse semigroup  topologies of semigroups of
partial homeomorphisms of some classes of topological spaces.  In
this context the results of our paper yield some notable results
about the topological behavior of the finite rank symmetric
inverse semigroups sitting inside larger function space
semigroups, or larger semigroups in general.  For example, under
reasonably general conditions, the inverse semigroup of partial
finite bijections $\mathscr{I}_\lambda^n$ of rank $\leqslant n$ is
a closed subsemigroup of a topological semigroup which contains
$\mathscr{I}_\lambda^n$ as a subsemigroup.

The class of semigroups that we consider in this paper is a
general class of semigroups that is modeled on the semigroup
$\J_\lambda^n$ and includes it as a special case. In Section 2 we
take a closer look at this class including some categorical
properties and a more general example than $\J_\lambda^n$.

A question of interest over the years has been to identify classes
of semigroups that can be embedded in compact semigroups and
classes that resist such embeddings. In Section 3 we consider this
question in the context of the class of semigroups we are
considering.

In this paper all topological spaces will be assumed to be
Hausdorff. We shall follow the terminology of~\cite{CHK, CP,
Engelking1989, Petrich1984, Ruppert1984}. If $A$ is a subset of a
topological space $X$, then we denote the closure of the set $A$
in $X$ by $\operatorname{cl}_X(A)$, or simply $\overline A$ if $X$
is obvious from context. By
\begin{equation*}
\left(%
\begin{array}{cccc}
  x_1 & x_2 & \cdots & x_n \\
  y_1 & y_2 & \cdots & y_n \\
\end{array}%
\right)
\end{equation*}
we denote a partial one-to-one transformation which maps $x_1$
onto $y_1$, $x_2$ onto $y_2$, $\ldots$, and $x_n$ onto $y_n$, and
by $0$ the empty transformation.

\section{Semigroup closures}
\begin{definition} \label{D:unstable} A subset $D$ of a semigroup
$S$ is said to be $\omega$-\emph{unstable} if $D$ is infinite and
for any $a\in D$ and infinite subset $B\subseteq D$, we have
$aB\cup Ba\not\subseteq D$.
\end{definition}

The next lemma gives a basic example of $\omega$-unstable sets.
\begin{lemma}\label{L1} For $\lambda$ infinite,
$D:=\J_\lambda^n\setminus {\J}_\lambda^{n-1}$ is an
$\omega$-unstable subset of  ${\J}_\lambda^n $. Indeed for $a\in
D$ and $B\subseteq D$ of cardinality at least $n!+1$, $aB\cup
Ba\not\subseteq D$.
\end{lemma}

\begin{proof} In order for $aB$ to be contained in $D$, it must be
the case that the range of each member of $B$ equals the domain of
$a$ and in order for $Ba$ to be contained in $D$, it must be the
case the every member of $B$ have the same domain as the range
$a$.  Thus there are only $n!$ possibilities for members of $B$.
\end{proof}

Recall that a semigroup $S$ is a \emph{semitopological semigroup}
if it is equipped with a Hausdorff topology for which all left
translation maps $\lambda_s$ and all right translation maps
$\rho_s$ are continuous~\cite{Ruppert1984}.  In this case we say
equivalently that the multiplication is \emph{separately
continuous}.  We come now to the crucial lemma for all that
follows.

\begin{lemma} \label{L2}
Let $S$ be a semitopological semigroup, let $T$ be a subsemigroup,
let $I$ be an ideal of $T$, and assume $D:=T\setminus I$ is
$\omega$-unstable.  If $s\in S$ is a limit point of $T$, then for
each $t\in\overline T$, either $st\in \overline I$ or $ts\in
\overline I$.
\end{lemma}

\begin{proof}  We have $s\in \overline{T}=\overline{D}\cup
\overline{I}$. If $s\in\overline{I}$, then by separate continuity
of multiplication for each $t\in T$, $st\in\overline{I}t\subseteq
\overline{It}\subseteq \overline{I}$.  Thus
$sT\subseteq\overline{I}$, and by continuity of left translation
by $s$, $s\overline{T}\subseteq \overline{I}$.  Similarly
$\overline{T}s\subseteq \overline{I}$.

For the case that $s\in\overline D$, but $s\notin \overline I$,
suppose for some $t\in T$ that $st, ts\in W:=S\setminus \overline
I$.  Then $t\notin I$, for otherwise $st, ts\in\overline I$, and
therefore $t\in D$. Using the continuity of left and right
translation by $t$, we find an open set $U$ containing $s$ such
that $Ut\cup tU\subseteq W$.  Since $s$ is a limit point of $T$
and $S$ is Hausdorff, the set $B:=T\cap (U\setminus \overline I)$
is infinite.  Since $B\subseteq T\setminus \overline I\subseteq D$
and $D$ is $\omega$-unstable, either $Bt$ or $tB$ meets $I$.  But
$tB\cup Bt\subseteq tU\cup Ut\subseteq W$, and $W$ misses $I$, a
contradiction. We conclude that for all $t\in T$, either
$st\in\overline I$ or $ts\in\overline I$. Thus the closed set
$(\lambda_s)^{-1}(\overline I)\cup(\rho_s)^{-1}(\overline I)$
contains $T$ and hence $\overline T$, which completes the proof.
\end{proof}

\begin{corollary} \label{C3}
Let $S$ be a semitopological semigroup, let $T$ be a subsemigroup,
let $I$ be an ideal of $T$, and assume $D:=T\setminus I$ is
$\omega$-unstable.  If $x,y\in \overline T$ and $xyx=x$, then
either $x\in \overline I$ or $x\in T$ and $x$ is an isolated point
of $\overline T$.
\end{corollary}

\begin{proof}
Suppose that $x\notin \overline I$. If $x\notin T$, then $x$ must
be a limit point of $T$.  By Lemma \ref{L2} either $xy\in
\overline I$ or $yx\in\overline I$.  But then
$x=(xy)x=x(yx)\in\overline I$, since $\overline I$ is an ideal of
$\overline T$  by separate continuity.  This contradicts our
assumption that $x\notin\overline I$. Thus $x\in
T\setminus\overline I$.  If $x$ is a limit point of $T$, then
Lemma \ref{L2} would again imply that $x=(xy)x=x(yx)\in \overline
I$, a contradiction.  Thus $x$ must be an isolated point of $T$
and hence also of $\overline T$.
\end{proof}

\begin{definition}
An \emph{ideal series} (see, for example, \cite{CP}) for a
semigroup $S$ is a chain of ideals
\begin{equation*}
 %$$
I_0\subseteq I_1\subseteq I_2\subseteq \ldots \subseteq I_m=S.
 %$$
\end{equation*}
We  call the ideal series \emph{tight} if $I_0$ is a finite set
and $D_k:=I_k\setminus I_{k-1}$ is an $\omega$-unstable subset for
each $k=1,\ldots, m$.
\end{definition}

\begin{example}\label{E1}
It follows from Lemma \ref{L1} that  for an infinite cardinal
$\lambda$, $\{\emptyset\}\subseteq
\J_\lambda^1\subseteq\ldots\subseteq \J_\lambda^m$ is a tight
ideal series for $S:=\J_\lambda^m$.
\end{example}

Recall that an element $x$ of a semigroup $S$ is regular if there
exists $y\in S$ such that $xyx=x$ and that $S$ is regular if every
element is regular~\cite{CP}.  If $xyx=x$, then it is
straightforward to verify that $x'=yxy$ is an inverse for $x$,
i.e., $xx'x=x$ and $x'xx'=x'$.  If $x$ belongs to an ideal $I$,
then $x'=yxy$ also belongs to the ideal, and it follows that an
ideal of a regular semigroup is regular.

\begin{proposition}\label{P4}
Let $S$ be a semitopological regular semigroup that admits a tight
ideal series $I_0\subseteq \ldots,\subseteq I_m=S$. Then each
$I_k$ is closed in $S$ and each member of $S\setminus I_{m-1}$ is
an isolated point of $S$.
\end{proposition}

\begin{proof}
We first prove by finite induction that $I_k$ is closed in $S$ for
each $k$. First note that $I_0$ is closed since it is finite and
$S$ is Hausdorff.  Assume that $I_{k-1}$ is closed for some $k$,
$1\leq k\leq m$.  If $s\in\overline{I_k}$, then by regularity of
the ideal $\overline{I_k}$ and Corollary \ref{C3} $s\in
\overline{I_{k-1}}=I_{k-1}$ or $s\in I_k$.  In either case $s\in
I_k$, so $I_k$ is closed.  By induction $I_k$ is closed for all
$k$, $0\leq k\leq m$.  The last assertion now follows from
Corollary \ref{C3}.
\end{proof}

By Example \ref{E1} and Proposition \ref{P4} we have the following

\begin{corollary}\label{C4.5}
Let $\lambda\geqslant\omega$ and let $n$ be any positive integer.
If $\tau$ is a topology on $\mathscr{I}_\lambda^{n}$ such that
$(\mathscr{I}_\lambda^{n},\tau)$ is a semitopological semigroup,
then every element
$\alpha\in\mathscr{I}_\lambda^{n}\setminus\mathscr{I}_\lambda^{n-1}$
is an isolated point of the topological space
$(\mathscr{I}_\lambda^{n},\tau)$.
\end{corollary}

The next proposition establishes that semigroups with tight ideal
series are restricted in regard to the type of semigroups in which
they may be densely embedded.

\begin{proposition}\label{P5}
Let $S$ be a semitopological semigroup, let $T$ be a subsemigroup,
and let $\overline T$ be its closure in $S$.  If $T$ admits a
tight ideal series, then any regular element of $\overline T$, in
particular any idempotent, must already be in $T$.
\end{proposition}

\begin{proof}
Let $x$ be a regular element of $\overline T$, and let
$I_0\subseteq \ldots\subseteq I_m=T$ be a tight ideal series for
$T$. Let $k$ be the smallest index such that $x\in
\overline{I_k}$.  We claim $x\in I_k$.  If $k=0$, then we are done
since $I_0$ is finite, hence closed.  So we assume $k\geq 1$. By
hypothesis there exists $y\in\overline T$ such that $xyx=x$. Then
$x':=yxy$ satisfies $xx'x=x$ and $x'=yxy\in \overline{I_k}$, since
$\overline{I_k}$ is an ideal of $\overline T$.  We apply Corollary
\ref{C3}  (with $T=I_k$ and $I=I_{k-1}$) and conclude that $x\in
I_k\subseteq T$.
\end{proof}

Recall that a  semigroup is an \emph{inverse semigroup} if it is a
regular semigroup in which each element $x$ has a unique inverse
$x'$~\cite{Petrich1984}.  The semigroups $\J_\lambda$ and
$\mathscr{I}_\lambda^n$ are inverse semigroups.

\begin{proposition} \label{P6}
Let $S$ be a semitopological inverse semigroup for which the
inversion map $x\mapsto x'$ is continuous.  If $T$ is an inverse
subsemigroup that admits a tight ideal series, then $T$ is closed
in $S$.
\end{proposition}

\begin{proof}
By separate continuity $\overline T$ is a subsemigroup.  Since $T$
is closed under the inversion map and the inversion map is
continuous, one readily sees that $\overline T$ is closed under
inversion, i.e., $\overline T$ is an inverse, hence regular,
subsemigroup.   Proposition \ref{P5} then yields that $\overline
T\subseteq T$, i.e., $T$ is closed.
\end{proof}

Proposition \ref{P6} applies directly to the symmetric inverse
semigroup $\J_\lambda^n$ for $\lambda$ infinite and $n$ a positive
integer and yields the following corollary.

\begin{corollary} \label{C7}
Let  $S$ be a semitopological inverse semigroup for which the
inversion map $x\mapsto x'$ is continuous.  If  $($an isomorphic
copy of$)\ \J_\lambda^n$ is a subsemigroup of $S$, then it is a
closed subset of $S$.
\end{corollary}

\begin{definition}[\cite{GutikPavlyk2001, Stepp1969}]\label{def1}
Let $\mathfrak{S}$ be a class of topological semigroups.
A~topological semigroup $S\in\mathfrak{S}$ is called
\emph{$H$-closed in the class $\mathfrak{S}$} if $S$ is a closed
subsemigroup of any topological semigroup $T\in{\mathfrak{S}}$
which contains $S$ as a subsemigroup. If $\mathfrak{S}$ coincides
with the class of all topological semigroups, then the semigroup
$S$ is called \emph{$H$-closed}.
\end{definition}

We remark that in~\cite{Stepp1969} the $H$-closed semigroups are
called {\it maximal}.

\begin{definition}[\cite{GutikPavlyk2001, Stepp1975}]\label{def2}
Let $\mathfrak{S}$ be a class of topological semigroups.
A~semigroup $S$ is called \emph{algebraically closed in the class}
$\mathfrak{S}$ if for any topology $\tau$ on $S$ such that
$(S,\tau)\in\mathfrak{S}$ we have that $(S,\tau)$ is an $H$-closed
topological semigroup in the class $\mathfrak{S}$. If
$\mathfrak{S}$ coincides with the class of all topological
semigroups, then the semigroup $S$ is called \emph{ algebraically
closed}.
\end{definition}

We have immediately from Corollary \ref{C7} the following
corollaries.

\begin{corollary}\label{C8}
For any infinite cardinal $\lambda$ and positive integer $n$, the
semigroup $\J_\lambda^n$ is algebraically closed in the class of
topological inverse semigroups $($inverse semigroups that are
topological semigroups with continuous inversion$)$.
\end{corollary}

\begin{corollary}\label{C9}
Let $n$ be any positive integer and let $\tau$ be any inverse
semigroup topology on $\mathscr{I}_\lambda^{n}$. Then
$(\mathscr{I}_\lambda^{n},\tau)$ is an $H$-closed topological
inverse semigroup in the class of topological inverse semigroups.
\end{corollary}

The following example implies that for all
$\lambda\geqslant\omega$, the semigroup $\mathscr{I}_\lambda^{k}$
with the discrete topology is not $H$-closed in the class of all
locally compact topological semigroups, for any positive integer
$k$.

\begin{example}\label{example5} We fix any positive integer $k$.
Let $a\notin \mathscr{I}_\omega^{k}$. Let
$S=\mathscr{I}_\omega^{k}\cup\{ a\}$. We put
\begin{equation*}
a\cdot a=a\cdot x=x\cdot a=0
\end{equation*}
for all $x\in \mathscr{I}_\omega^{k}$.

We further enumerate the elements of the set $\omega$ by natural
numbers. Let
\begin{equation*}
 A_m=\left\{\left(%
\begin{array}{c}
  2l-1 \\
  2l \\
\end{array}%
\right)\;\Big| \;l\geqslant m\right\}
\end{equation*} for each positive integer $m$. A
topology $\tau$ on $S$ is now defined as follows:
\begin{enumerate}
    \item[1)] all points of $\mathscr{I}_\omega^{k}$ are isolated in
    $S$; and
    \item[2)] $\mathscr{B}(a)=\{ U_n(a)=\{ a\}\cup A_n\mid
    n=1,2,3,\ldots\}$ is the base of the topology $\tau$ at the
    point $a\in S$.
\end{enumerate}

Then
\begin{enumerate}
\item[a)]
for all $\left(%
\begin{array}{cccc}
  x_1 & x_2 & \cdots & x_i \\
  y_1 & y_2 & \cdots & y_i \\
\end{array}%
\right)\in \mathscr{I}_\omega^{k}$ and $n\geqslant\max\{ x_1, x_2,
\ldots, x_i, y_1, y_2, \ldots, y_i\}$ we have
\begin{equation*}
\left(%
\begin{array}{cccc}
  x_1 & x_2 & \cdots & x_i \\
  y_1 & y_2 & \cdots & y_i \\
\end{array}%
\right)
 \cdot U_n(a)= U_n(a)\cdot
\left(%
\begin{array}{cccc}
  x_1 & x_2 & \cdots & x_i \\
  y_1 & y_2 & \cdots & y_i \\
\end{array}%
\right)=\{ 0\};
\end{equation*}

\item[b)] $U_n(a)\cdot U_n(a)=U_n(a)\cdot\{ 0\}=\{ 0\}\cdot
U_n(a)=\{ 0\}$ for any positive integer $n$; and

\item[c)] $U_n(a)$ is a compact subset of $S$ for each positive
integer $n$.
\end{enumerate}

Therefore $(S,\tau)$ is a locally compact topological semigroup.
Obviously $\mathscr{I}_\omega^{k}$ is not a closed subset of $(S,
\tau)$.
\end{example}

The following example shows that for all $\lambda\geqslant\omega$,
the semigroup $\mathscr{I}_\lambda^{\infty}:=\bigcup_n
\J_\lambda^n$ with the discrete topology is not $H$-closed in the
class of all topological inverse semigroups.

\begin{example}\label{example6}
Let $\lambda\geqslant\omega$ and let $\tau_d$ be the discrete
topology on the semigroup $\mathscr{I}_\lambda^{\infty}$.

For any $\varepsilon\in E(\mathscr{I}_\lambda^{\infty})$ we define
\begin{equation*}
M(\varepsilon)=\{\chi\in\mathscr{I}_\lambda^{\infty}\mid
\varepsilon\chi=\chi\varepsilon=\varepsilon\}.
\end{equation*}

Let $S$ be the semigroup $\mathscr{I}_\lambda^{\infty}$ with the
adjoined identity $\iota$. We now define a topology $\tau_S$ on
the semigroup $S$ as follows:
\begin{itemize}
    \item[$(i)$] $\chi$ is an isolated point in $S$ for all
    $\chi\in\mathscr{I}_\lambda^{\infty}$; and
    \item[$(ii)$] the family
    \begin{equation*}
      \mathscr{B}(\iota)=\left\{ U_\varepsilon(\iota)=\{\iota\}\cup
      M(\varepsilon)\mid \varepsilon\in
      E(\mathscr{I}_\lambda^{\infty})\right\}
    \end{equation*}
    is the base of the topology $\tau_S$ at the point $\iota$.
\end{itemize}

The definition of the family $\mathscr{B}(\iota)$ implies that
$\iota$ is not an isolated point of a topological space
$(S,\tau_S)$ and the restriction of the topology $\tau_S$ on the
set $\mathscr{I}_\lambda^{\infty}$ coincides with the topology
$\tau_d$.

Obviously, this is sufficient to show that the semigroup operation
on $(S,\tau_S)$ is continuous in the following cases:
\begin{itemize}
    \item[$(i)$] $\iota\iota=\iota$; and
    \item[$(ii)$] $\iota\chi=\chi\iota=\chi$ for all
    $\chi\in\mathscr{I}_\lambda^{\infty}$.
\end{itemize}

In case $(i)$ we have
\begin{equation*}
    U_\varepsilon(\iota)\cdot U_\varepsilon(\iota)\subseteq
    U_\varepsilon(\iota).
\end{equation*}

In case $(ii)$ we denote
\begin{equation*}
  \chi=
\left(%
\begin{array}{cccc}
  x_1 & x_2 & \cdots & x_n \\
  y_1 & y_2 & \cdots & y_n \\
\end{array}%
\right).
\end{equation*} Then we put
\begin{equation*}
 K=\{ x_1\}\cup\{ x_2\}\cup\cdots\cup\{ x_n\}\cup\{
y_1\}\cup\{ y_2\}\cup\cdots\cup\{ y_n\}.
\end{equation*}
Let be $K=\{ a_1, a_2, \ldots, a_k\}$. Obviously $k\leqslant n$.
We define
\begin{equation*}
    \varepsilon=
\left(%
\begin{array}{cccc}
  a_1 & a_2 & \cdots & a_k \\
  a_1 & a_2 & \cdots & a_k \\
\end{array}%
\right).
\end{equation*}
Then we have $\chi\varepsilon=\varepsilon\chi=\chi$ and hence
\begin{equation*}
U_\varepsilon(\iota)\cdot\chi=
U_\varepsilon(\iota)\cdot\chi=\{\chi\}.
\end{equation*}

Since
$\left(U_\varepsilon(\iota)\right)^{-1}=U_\varepsilon(\iota)$, we
have that $(S,\tau_S)$ is a topological inverse semigroup which
contains $\mathscr{I}_\lambda^{\infty}$ as  dense inverse
subsemigroup.
\end{example}

\section{Semigroups with tight ideal series}

We have seen in the previous section that semigroups admitting a
tight ideal series have interesting closure properties in larger
semigroups.  In this section we take a brief closure look at this
class of semigroups, primarily to see that such semigroups extend
signficantly beyond the finite rank symmetric inverse semigroups.

\begin{lemma}\label{L2.1}
The class of semigroups admitting a tight ideal series is closed
under finite products.
\end{lemma}

\begin{proof}
It suffices to check for the case $n=2$.  Let $S$ have a tight
ideal series $I_0\subseteq \ldots \subseteq I_m=S$ and $T$ have a
tight ideal series $J_0\subseteq \ldots \subseteq J_n=T$.  Set
$K_i=I_i\times J_0$ for $0\leq i\leq m$ and $K_i=S\times J_{i-m}$
for $m<i\leq n+m$.  Then any infinite $B\subseteq K_{i+1}\setminus
K_i$ has an infinite projection into either $S$ for $i\leq m$ and
into $T$ for $i>m$, and since multiplication is coordinatewise, it
directly follows that $aB\cup Ba$ meets $K_i$.
\end{proof}

\begin{lemma}\label{L2.2}
Let $h\colon S\to T$ be a surjective semigroup homomorphism such
that each point inverse $h^{-1}(t)$ is finite.  If $S$ has a tight
ideal series, then so does $T$.
\end{lemma}

\begin{proof}
It is easy to see that $h^{-1}(I_0)\subseteq\ldots\subseteq
h^{-1}(I_m)=S$ is a tight ideal series for $S$ if
$I_0\subseteq\ldots\subseteq I_m$ is one for $T$.
\end{proof}

\begin{example}\label{2.3}
Let $\{X_i\mid i\in I\}$ be a collection of finite, pairwise
disjoint sets, and let $X=\bigcup_{i\in I} X_i$.  We consider the
semigroup $\mathscr{PT}(X,I)$ of partial functions on $X$ defined
in the following way.  First choose a subset $J\subset I$ and let
$\alpha\colon J\to K$ be a bijection to another subset of $I$. A
function $f\colon \bigcup_{i\in J} X_i\to \bigcup_{i\in K}$ is by
definition in $\mathscr{PT}(X,I)$ if and only if $x\in X_i$
implies $f(x)\in X_{\alpha(i)}$.  All such partial functions
ranging over all such $\alpha$ and corresponding $f$ form the
subsemigroup $\mathscr{PT}(X,I)$ under composition, a subsemigroup
of the semigroup of all partial functions on $X$. We define such a
function to have rank $n$ if the cardinality of the range (and
hence domain) of $\alpha$ has cardinality $n$ and
$\mathscr{PT}(X,I)^n$ to be all functions of rank less than or
equal to $n$.  We then have an ideal series for the semigroup
$\mathscr{PT}(X,I)^n$ defined by $I_k=\mathscr{PT}(X,I)^k$ for
$0\leq k\leq n$.

There is a surjective homomorphism $h$ from $\mathscr{PT}(X,I)^n$
to $\J_I^n$, which assigns to $f\in\mathscr{PT}(X,I)^n$ the
corresponding $\alpha\colon J\to K$ between the index sets.  One
checks directly that this homomorphism has finite point inverses,
and hence it follows from Lemma \ref{L2.2} that the ideal series
$I_k=\mathscr{PT}(X,I)^k$ for $0\leq k\leq n$ is tight.

We remark that the semigroup $\mathscr{PT}(X,I)^n$ is regular, but
not an inverse semigroup.  It is well known that the semigroup of
all transformations $\mathscr{T}(X)$ is regular \cite{CP}, and
essentially the same proof yields the regularity of
$\mathscr{PT}(X,I)^n$.  Thus the principal results of Section 1
may be applied to the semigroup $\mathscr{PT}(X,I)^n$.
\end{example}

\section{Compact embeddings}

Several authors have considered the problem of showing that
various specific semigroups or classes of semigroups do or do not
embed into compact semigroups. For example one can use the
Swelling Lemma to show that the bicyclic semigroup does not admit
an embedding into a compact topological semigroup \cite{CHK}.
Closer to our current investigations, it was shown by Gutik and
Pavlyk in \cite{GutikPavlyk2005} that an infinite topological
semigroup of $\lambda{\times}\lambda$-matrix units $B_{\lambda}$
does not embed into a compact topological semigroup, every
non-zero element of $B_{\lambda}$ is an isolated point of
$B_{\lambda}$, and $B_{\lambda}$ is algebraically closed in the
class of topological inverse semigroups.  (This is essentially a
special case of results of this paper for $\J_\lambda^1$.)
However, we add a new wrinkle to earlier investigations by showing
that certain partially compact embeddings do not exist, more
precisely that the closure of certain embedded $\D$-classes cannot
be compact.

Recall the Green's relations on a semigroup $S$.  Two elements are
$\mathcal{L}$-equivalent if they generate the same principal left
ideal, i.e., $s\mathcal{L} t$ if $\{s\}\cup Ss=\{t\}\cup St$, and
$\mathcal{R}$ related if they generate the same principal right
ideal.  (In the case that $s$ is regular the principal left ideal
reduces to $Ss$ since $s=ss's\in sS$.)  The join of the
equivalence relations $\mathcal{L}$ and $\mathcal{R}$ is denoted
$\D$. It is a standard semigroup result that $\D$ is alternatively
given by the relational compositions $\mathcal{L}\circ\mathcal{R}=
\mathcal{R}\circ\mathcal{L}$ \cite[Section 2.1]{CP}.   A
$\D$-equivalence class $D$ is called a \emph{regular $\D$-class}
if it contains a regular element.  This is the case if and only if
each $\mathcal{L}$-class and each $\mathcal{R}$-class contained in
$D$ contains at least one idempotent if and only if every element
of $D$ is regular  \cite[Chapter 2.3]{CP}. Furthermore, each
inverse of a member of $D$ is back in $D$ \cite[Chapter 2.3]{CP}

We recall a useful fact about regular $\D$-classes.

\begin{lemma}\label{3.1}
Let $a,c\in D$, a regular $\D$-class in a semigroup $S$.  Then
there exist $s,t\in D$ such that $c=sat$.
\end{lemma}

\begin{proof}
Since $\D=\mathcal{R}\circ\mathcal{L}$, we may pick $b\in D$ such
that $a\mathcal{R}b$ and $b\mathcal{L} c$.  Pick an idempotent $e$
in the $\mathcal{R}$-class of $a$ and $u\in S$ such that $au=e$.
Since $eS=aS=bS$, we have also $ea=a$ and $eb=b$.  For $t=ub$, we
have $at=aub=eb=b$.  Furthermore, $t=ub\in Sb$ and $b=at\in St$,
so $t\mathcal{L}b$, and thus $t\in D$.  In a similar fashion one
finds $s\in D$ such that $c=sb$.  Then $c=sb=sat$ and $s,t\in D$.
\end{proof}

The next theorem is our main one on the non-existence of compact
embeddings of certain $\D$-classes. \begin{theorem}\label{3.1} Let
$S$ be a topological semigroup and let $T$ be a subsemigroup
having a tight ideal series $I_0\subseteq \ldots\subseteq I_m$.
If $D:=I_{k+1}\setminus I_k$ is a regular $\D$-class, then
$\overline D=\operatorname{cl}_S(D)$ is not compact.
\end{theorem}

\begin{proof}
Suppose the contrary.  Then the infinite set $D$ has a limit point
$x$ in the compact set $\overline D$.  Let $x_\alpha$ denote a net
in $D\setminus\{x\}$ converging to to $x$.  For each $\alpha$,
pick an inverse $x_\alpha'$, which must again be in $D$.  By
compactness of $\overline D$, some subnet of $x_\alpha'$ (which we
again label $x_\alpha'$) must converge to some $y\in \overline D$.
By continuity of multiplication $x_\alpha=x_\alpha
x_\alpha'x_\alpha\to xyx$, and by uniqueness of limits $x=xyx$.
Thus $x$ is a regular element.

Fix some $a\in D$.  By Lemma \ref{3.1} for each $\alpha$, there
exists $s_\alpha,t_\alpha\in D$ such that $s_\alpha x_\alpha
t_\alpha=a$.  Again passing to convergent subnets, we have
$s_\alpha\to s\in\overline D$, $t_\alpha\to t\in\overline D$, and
$a=s_\alpha x_\alpha t_\alpha\to sxt$.  Therefore $a=sxt$.   In a
similar fashion, we can write $x_\alpha=u_\alpha a v_\alpha$ and
conclude $x=uav$ for $u,v\in\overline D$.

It follows from Corollary \ref{C3} that $x\in \overline{I_k}$.
Let $j$ be the smallest index such that $x\in\overline{I_j}$.
Then $j\ne 0$, for otherwise since $I_0$ is finite, hence closed,
$x\in I_0$, and hence $a=sxt\in I_0$, contradicting the fact that
$a\notin I_k$, which contains $I_0$.  Thus $j\geq 1$.  If $x\notin
I_j$, then $x$ would be a limit point of $I_j$, and hence a limit
point of $I_j\setminus I_{j-1}$, since
$x\notin\overline{I}_{j-1}$.  Again by Corollary \ref{C3} it would
follow that $x\in \overline{I}_{j-1}$, a contradiction.  We
conclude that $x\in I_j$, and then that $a=sxt\in \overline{I_j}$,
since $\overline{I_j}$ is an ideal of $\overline T$.  Applying
Corollary \ref{C3} to the regular element $a$, we conclude either
that $a\in I_j$, an impossibility since $I_j\subseteq I_k$, or
$a\in\overline{I}_{j-1}$.  The latter would imply $x=uav\in
\overline{I}_{j-1}$, which we have just seen is not the case.
Thus we have reached a contradiction to our assumption that
$\overline D$ is compact.
\end{proof}

Since it is well-known and direct to verify that the sets
$D=\J_\lambda^k\setminus \J_\lambda^{k-1}$ are $\D$-classes for
$1\leq k\leq n$ in the semigroup $\J_\lambda^n$, we have the
following corollary, which generalizes the previously mentioned
result of Gutik and Pavlyk in \cite{GutikPavlyk2005}.

\begin{corollary}
For an infinite cardinal $\lambda$ and positive integer $n$, if
$\J_\lambda^n$ is a subsemigroup of a topological semigroup $S$,
it cannot be the case that that
$\operatorname{cl}_S(\J_\lambda^k\setminus \J_\lambda^{k-1})$ is
compact for $1\leq k\leq n$.
\end{corollary}

We close with a theorem on $\J_\lambda^\infty$.

\begin{theorem}\label{theorem12}
For any infinite cardinal $\lambda$ there exists no topology
$\tau$ on $\mathscr{I}_\lambda^{\infty}$ such that
$(\mathscr{I}_\lambda^{\infty},\tau)$ is a compact semitopological
semigroup.
\end{theorem}

\begin{proof}
Suppose to the contrary, that there exists a topology $\tau$ on
$\mathscr{I}_\lambda^{\infty}$ such that
$(\mathscr{I}_\lambda^{\infty},\tau)$ is a compact semitopological
semigroup.

The definition of the semigroup $\mathscr{I}_\lambda^{\infty}$
implies that for any idempotent
$\varepsilon\in\mathscr{I}_\lambda^{\infty}$ there exists an
idempotent $\phi\in\mathscr{I}_\lambda^{\infty}$ such that
$\varepsilon\phi=\phi\varepsilon=\varepsilon$ and
$\phi\neq\varepsilon$, i.~e. $\varepsilon<\phi$. Therefore there
exists a subsets of idempotents $A=\{\varepsilon_1,
\varepsilon_2,\ldots, \varepsilon_n,\ldots\}$ in
$\phi\in\mathscr{I}_\lambda^{\infty}$ such that
\begin{equation*}
0<\varepsilon_1< \varepsilon_2<\ldots<\varepsilon_n<\ldots .
\end{equation*}
Without loos of generality we can assume that
$\varepsilon_k\in\mathscr{I}_\lambda^{k}\setminus
\mathscr{I}_\lambda^{k-1}$, for any $k=2,3,4,\ldots$ and
$\varepsilon_1\in\mathscr{I}_\lambda^{1}\setminus\{ 0\}$. Then
Corollary~\ref{C4.5} implies that the idempotent $\varepsilon_k$
has an open neighbourhood $U(\varepsilon_k)$ such that
$U(\varepsilon_k)\cap\mathscr{I}_\lambda^{k}=\{\varepsilon_k\}$
for all $k=1,2,3,\ldots$. Since the translations in
$(\mathscr{I}_\lambda^{n},\tau)$ are continuous maps, the set
\begin{equation*}
U_{l}(\varepsilon_k)=\{\beta\in\mathscr{I}_\lambda^{\infty}\mid
\beta\varepsilon_k=\varepsilon_k\}
\end{equation*}
is clopen in the topological space
$(\mathscr{I}_\lambda^{\infty},\tau)$ for all $k=1,2,3,\ldots$. We
define the family
\begin{equation*}
\mathscr{O}=\{O_k\mid k=1,2,3,\ldots\}
\end{equation*}
as follows:
\begin{itemize}
    \item[$(i)$] $O_1=\mathscr{I}_\lambda^{\infty}\setminus
    U_{l}(\varepsilon_1)$; and
    \item[$(ii)$] $O_k=U_{l}(\varepsilon_{k-1})\setminus
    U_{l}(\varepsilon_{k})$ for all $k=2,3,4, \ldots$.
\end{itemize}
Obviously, the family $\mathscr{O}$ is a clopen cover of the
topological space $(\mathscr{I}_\lambda^{\infty},\tau)$, which
does not contain a finite subcover, a contradiction.
\end{proof}

%%%%%%%%%%%%%%%%%%%%%%%%%%%%%%%%%%%%%%%%%%%%%%%%%%%%%%%%%%%

\section*{Acknowledgements}

This research was supported by SRA grants
P1-0292-0101-04, J1-9643-0101 and 
BI-RU/08-09-002.

%%%%%%%%%%%%%%%%%%%%%%%%%%%%%%%%%%%%%%%%%%%%%%%%%%%%%%%%%%%%

\end{document}